\documentclass[11pt]{article}
\usepackage{amsmath}
\usepackage{amssymb}
\title{\bf{A Few Equivalences of Wall-Sun-Sun Prime Conjecture}}
\author{\bf{Arpan Saha}$^{\dag}$ \and \bf{Karthik C S}$^{\dag}$}

\date{}
\begin{document}
\maketitle
\begin{abstract}
\noindent In this paper, we prove a few lemmas concerning Fibonacci numbers modulo primes and provide a few statements that are equivalent to Wall-Sun-Sun Prime Conjecture. Further, we investigate the conjecture through heuristic arguments and propose a couple of additional conjectures for future research.
\end{abstract}
\let\thefootnote\relax\footnotetext{$\dag$ Arpan Saha (arpan[underscore]saha[at]iitb[dot]ac[dot]in) and Karthik C S (karthik[dot]c[dot]s[at]iitb[dot]ac[dot]in) are sophomore undergraduates at the Indian Institute of Technology (IIT), Bombay.}

\section{Introduction}
The Fibonacci sequence $\{F_n\}_{n\ge0}$ (defined as $F_0 = 0, F_1 = 1$ and $F_n = F_{n-1}+F_{n-2}$ for all $n\ge2$) has harbored great interest owing to its wide occurrence in combinatorial problems such as that regarding the number of ways of tiling a $2\times n$ rectangle with $2\times 1$ dominoes and other numerous properties it exhibits. For instance, $F_{m+n}=F_{m-1}F_n+F_mF_{n+1}$ and, as E.~Lucas had discovered, $F_{\gcd(m,n)}=\gcd(F_m,F_n)$. Moreover, the work of D.~D.~Wall, Z.~H.~Sun and Z.~W.~Sun \cite{Wall}\cite{Sun-Sun} regarding what has come to be known as the Wall-Sun-Sun Prime Conjecture, had demonstrated intimate links between the Fibonacci sequence and Fermat's Last Theorem \cite{Sun-Sun}. Though the latter i.e. Fermat's Last Theorem was proved in 1995 by Andrew Wiles and Richard Taylor \cite{FLT}\cite{Hecke}, the Wall-Sun-Sun Prime Conjecture continues to generate interest. This may be partly due to the fact that the Fibonacci sequence is interesting in its own right and partly due to the fact that it \emph{may} lead to a relatively elementary approach to Fermat's Last Theorem as compared to Wiles' proof involving bijections between elliptic and modular forms.
\\

The Wall-Sun-Sun Prime Conjecture is as follows:
\\

\emph{\bf{Statement 1:}} There does not exist a prime $p$ such that
$$
p^2\mid F_{p-\left(\frac{p}{5}\right)}
$$
where $\left(\frac{p}{5}\right)$ is the Legendre symbol i.e.

$$
\left(\frac{p}{5}\right) = \left\{
\begin{array}{rcll}
1 & \mbox{if} & p \equiv \pm 1 &\mbox{ (mod } 5)\\
-1 & \mbox{if} & p \equiv \pm 2 &\mbox{ (mod } 5)\\
0 & \mbox{if} & p \equiv 0 &\mbox{ (mod } 5)\\
\end{array}
\right.
$$

(Such primes shall henceforth be referred to as Wall-Sun-Sun primes.)

We shall provide a few statements equivalent to the above. But we will first require a few definitions and results.
\\

\emph{\bf{Definition:}} For a given positive integer $n$, $\kappa(n)$ is the least positive integer $m$ such that $n\mid F_m$.
\\

\emph{\bf{Definition:}} For a given positive integer $n$, $\pi(n)$ is the least positive integer $m$ such that $n\mid F_m$ and $F_{m+1}\equiv 1$ (mod $n$). This is often referred to as the Pisano period.
\\
\begin{sloppypar}
The existence of $\pi(n)$ for any positive integer $n$ follows from the Pigeonhole\ Principle and the well-ordering of positive integers \cite{Wall}; the existence of $\kappa(n)$ follows thence.
\end{sloppypar}
\section{Background Results}

We list here the results that we shall be using for demonstrating the equivalences discussed in the subsequent section.
\\

\emph{\bf{Lemma 1:}} Let $m$ and $n$ be positive integers. We claim that $n\mid F_m$ if and only if $\kappa(n)\mid m$.
\begin{sloppypar}
\emph{\bf{Proof:}} Both the necessity and sufficiency follow from the standard result due to E. Lucas that $\gcd(F_m,F_k) = F_{\gcd(m,k)}$ for any positive integers $m$ and $k$. Here, if $n\mid F_m$ then $n\mid \gcd(F_m,F_{\kappa(n)}) = F_{\gcd(m,\kappa(n))}$. But $F_{\kappa(n)}$, by definition is the least Fibonacci number divisible by $n$ (the Fibonacci numbers are an increasing sequence). So, $\kappa(n)\le\gcd(m,\kappa(n))$ which implies $\kappa(n)=\gcd(m,\kappa(n))$ i.e. $\kappa(n)\mid m$. Conversely, if $\kappa(n)\mid m$, then $\gcd(F_m,F_{\kappa(n)}) = F_{\gcd(m,\kappa(n))} = F_{\kappa(n)}$ i.e. $F_{\kappa(n)}\mid F_m$ from which we conclude $n\mid F_m$. $\blacksquare$
\\
\end{sloppypar}
\emph{\bf{Lemma 2:}} Let $l$ be the highest power of a positive integer $n$ dividing $F_{\kappa(n)}$. If $n\mid F_m$ for some positive integer $m$, then $n^l\mid F_m$.

\emph{\bf{Proof:}} As, $\gcd(F_m,F_{\kappa(n)}) = F_{\gcd(m,\kappa(n))} = F_{\kappa(n)}$ whenever $n$ divides $F_m$, we have $F_{\kappa(n)}\mid F_m$. Since $n^l\mid F_{\kappa(n)}$, we get $n^l\mid F_m$ as well. $\blacksquare$
\\

\emph{\bf{Lemma 3:}} Let $p$ be a prime. If $p\mid F_{pm}$ and $p\nmid F_p$, then $p\mid F_m$.

\emph{\bf{Proof:}} Let us recall Siebeck's formula for $F_{mn}$
\begin{equation*}
F_{mn} = \sum_{j=0}^{n}\binom{n}{j}F_jF_m^jF_{m-1}^{n-j}
\end{equation*}
where $m$ and $n$ are positive integers. Let us put $n=p$. We would then have
\begin{eqnarray*}
F_{mp} &=& \sum_{j=0}^{p}\binom{p}{j}F_jF_m^jF_{m-1}^{p-j}\\
&=& \sum_{j=1}^{p-1}\binom{p}{j}F_jF_m^jF_{m-1}^{p-j}+F_pF_m^p
\end{eqnarray*}
On taking the equality modulo p, we have
\begin{equation*}
F_{mp}\equiv F_pF_m^p \indent\mbox{ (mod } p)
\end{equation*}
as $p$ divides $\binom{p}{j}$
for $1\le j \le p-1$.

Now from Fermat's Little Theorem we obtain,
\begin{equation*}
F_{mp}\equiv F_pF_m \indent\mbox{ (mod } p)
\end{equation*}
from which the result follows. $\blacksquare$
\\

\emph{\bf{Lemma 4:}} Let $l$ be the highest power of a prime $p$ dividing $F_{\kappa(p)}$. We have,
 $$
 \kappa(p^{l+1}) = p\kappa(p)
 $$

\emph{\bf{Proof:}} We put $n=p$ and $m = \kappa(p)$ in Siebeck's formula.
\begin{equation*}
F_{p\kappa(p)} = \sum_{j=0}^{p}\binom{p}{j}F_jF_{\kappa(p)}^jF_{\kappa(p)-1}^{p-j}
\end{equation*}
Let $F_{\kappa(p)} = p^l\alpha$ for some positive integer $\gamma$ not divisible by $p$. Then,
\begin{eqnarray*}
F_{p\kappa(p)} &=& \sum_{j=0}^{p}\binom{p}{j}F_jp^{lj}\gamma^jF_{\kappa(p)-1}^{p-j}\\
&=& p^{l+1}\gamma F_{\kappa(p)-1}^{p-1}+p^{2l}\sum_{j=2}^{p}\binom{p}{j}F_jp^{l(j-2)}\gamma^jF_{\kappa(p)-1}^{p-j}
\end{eqnarray*}
As $p$ does not divide $\gamma$ or $F_{\kappa(p)-1}$ and we know that $l\ge1$, this implies that the highest power of $p$ dividing $F_{p\kappa(p)}$ is $l+1$.

Now, by Lemma 1, we gather that $\kappa(p^{l+1})\mid p\kappa(p)$. Let $a$ be a positive integer such that $p\kappa(p) = a\kappa(p^{l+1})$. Since $p^{l+1}\mid F_{\kappa(p^{l+1})}$, we have $p\mid F_{\kappa(p^{l+1})}$ and again by Lemma 1, we have $\kappa(p)\mid\kappa(p^{l+1})$. Let $b$ be a positive integer such that $\kappa(p^{l+1}) = b\kappa(p)$. So, $p\kappa(p)=ab\kappa(p)$ which is to say, $p = ab$. Now $a$ can be $p$ or $1$ as $p$ is prime. If $a = p$, then we would have $\kappa(p) = \kappa(p^{l+1})$ i.e. $p^{l+1}\mid F_{\kappa(p)}$. This contradicts our assumption, so we conclude $a = 1$ i.e.
$$
\kappa(p^{l+1}) = p\kappa(p) \indent \blacksquare
$$
\\

\emph{\bf{Lemma 5:}} For all positive integers $m$ and $n$,
\begin{equation*}
F_{m\kappa(n) + 1}\equiv F_{\kappa(n) + 1}^m \indent\mbox{ (mod } n^2)
\end{equation*}
\indent\emph{\bf{Proof:}} We proceed by induction. The base case $m = 1$ is trivial. For the rest, we shall invoke the standard result,
\begin{equation*}
F_{r+s} = F_rF_{s+1} + F_{r-1}F_s,~~~~~\forall r,s \in \mathbb{N}
\end{equation*}
Let us assume that the congruence holds for all positive integers $m~<~\lambda,~~ \lambda\in~\mathbb{N}$. Then,
\begin{eqnarray*}
F_{\lambda\kappa(n)+1} &=& F_{(\lambda-1)\kappa(n)+1+\kappa(n)}\\
&=& F_{(\lambda-1)\kappa(n)+1}F_{\kappa(n)+1}+F_{(\lambda-1)\kappa(n)}F_{\kappa(n)}\\
&\equiv& F_{\lambda\kappa(n)+1}^{\lambda-1}F_{\lambda\kappa(n)+1} \indent\mbox{ (mod } n^2) \indent\mbox{ (by inductive hypothesis) }\\
&\equiv& F_{\lambda\kappa(n)+1}^\lambda \indent\mbox{ ~~~~~~~~~~~(mod } n^2)
\end{eqnarray*}
Hence, the congruence holds for $m=\lambda$ as well. The Lemma is thus proved.~$\blacksquare$
\\

\emph{\bf{Lemma 6:}} Let $n$ be a positive integer and let $\Omega_n(z)$ denote the order of a positive integer $z$ modulo $n$. We have,
\begin{equation*}
\pi(n)=\kappa(n)\Omega_n(F_{\kappa(n)+1})
\end{equation*}
\indent\emph{\bf{Proof:}} The least positive integer $m$ for which,
\begin{eqnarray*}
F_{m\kappa(n) + 1}&\equiv& 1\indent\mbox{ (mod } n)\\
\mbox{i.e. }\indent F_{\kappa(n)+1}^m&\equiv& 1 \indent\mbox{ (mod } n)\indent\mbox{ (by Lemma 5)}
\end{eqnarray*}
is $\Omega_n(F_{\kappa(n)+1})$. Since $\gcd(F_{\kappa(n)+1}, F_{\kappa(n)})=1$, we have that $n$, which divides $F_{\kappa(n)}$, is relatively prime to $F_{\kappa(n)+1}$. Hence, $\Omega_n(F_{\kappa(n)+1})$ is well-defined. Now, from our definition of $\pi(n)$, the Lemma immediately follows. $\blacksquare$
\\

\emph{\bf{Lemma 7:}} Let $r$ and $n$ be positive integers and $p$ be a prime. Let $\alpha_r$ and $\beta_r$ be residues modulo $p$ such that
\begin{eqnarray*}
F_{r\pi(n)} &\equiv& \alpha_rp \indent\mbox{ ~~~~~(mod } p^2)\\
F_{r\pi(n)+1} &\equiv& \beta_rp+1 \indent\mbox{ (mod } p^2)
\end{eqnarray*}
We claim that
\begin{eqnarray*}
\alpha_r &\equiv& r\alpha_1 \indent\mbox{ (mod } p)\\
\beta_r &\equiv& r\beta_1 \indent\mbox{ (mod } p)
\end{eqnarray*}
\indent\emph{\bf{Proof:}} We first note that the above notation is well-defined as,
\begin{eqnarray*}
F_{r\pi(n)} &\equiv& 0 \indent\mbox{ (mod } p)\\
F_{r\pi(n)+1} &\equiv& 1 \indent\mbox{ (mod } p)
\end{eqnarray*}
We now proceed by induction. The base case $r=1$ is trivial. For the rest, we first assume that the Lemma holds for $r<\varrho$. From induction hypothesis, we have,
\begin{eqnarray*}
F_{(\varrho-1)\pi(p)} &\equiv& \alpha_{\varrho-1}p\\
&\equiv& (\varrho-1)\alpha_1p \indent\mbox{~~~~~ (mod } p^2)\\
F_{(\varrho-1)\pi(p)+1} &\equiv& \beta_{\varrho-1}p+1\\
&\equiv& (\varrho-1)\beta_1p+1 \indent\mbox{ (mod } p^2)
\end{eqnarray*}
Now for the inductive steps:
\begin{eqnarray*}
F_{\varrho\pi(p)} &\equiv& F_{(\varrho-1)\pi(p)+\pi(p)}\\
&\equiv& F_{\pi(p)-1}F_{(\varrho-1)\pi(p)}+F_{\pi(p)}F_{(\varrho-1)\pi(p)+1}\\
&\equiv& ((\beta_1-\alpha_1)p+1)(\varrho-1)\alpha_1p+\alpha_1p((\varrho-1)\beta_1p+1)\\
&\equiv& \alpha_1p\varrho \indent\mbox{ (mod } p^2)
\end{eqnarray*}
This gives us
\begin{equation*}
\alpha_{\varrho} \equiv \alpha_1\varrho \indent\mbox{ (mod } p)
\end{equation*}
Similarly,
\begin{eqnarray*}
F_{\varrho\pi(p)+1} &\equiv& F_{(\varrho-1)\pi(p)+1+\pi(p)}\\
&\equiv& F_{\pi(p)}F_{(\varrho-1)\pi(p)}+F_{\pi(p)+1}F_{(\varrho-1)\pi(p)+1}\\
&\equiv& (\varrho-1)\alpha_1^2p^2+(\beta_1p+1)((\varrho-1)\beta_1p+1)\\
&\equiv& \beta_1p\varrho+1 \indent\mbox{ (mod } p^2)
\end{eqnarray*}
This gives us
\begin{equation*}
\beta_{\varrho} \equiv \beta_1\varrho \indent\mbox{ (mod } p)
\end{equation*}
Thus the Lemma holds for all positive integers $r$. $\blacksquare$
\\

\emph{\bf{Lemma 8:}} Let $p$ be a prime. We claim that $\pi(p^2)$ equals either $\pi(p)$ or $p\pi(p)$.$^{\dag}$
\let\thefootnote\relax\footnotetext{$\dag$ Although this Lemma is very well established, we hope that we have given a new and more elementary proof.}

\emph{\bf{Proof:}} Firstly, we know that $\pi(n)\mid \pi(n^2)$ from Theorem 5 stated in \cite{Wall} and simple counting arguments. Hence, $\pi(p^2)$ is of the form $\xi\pi(p)$ where $\xi$ is a positive integer. We continue with the notation $\alpha_1$ and $\beta_1$ as introduced in Lemma 7 and investigate three cases:

\emph{Case (i):} \indent$\alpha_1 \not\equiv 0 $ \indent(mod $p$)

In this case, $\alpha_1$, $2\alpha_1$, $3\alpha_1$, ... $(p-1)\alpha_1$ are not congruent to zero modulo $p$. Thus, by Lemma 7, we have for all positive integers $\xi<p$,
\begin{equation*}
F_{\xi\pi(p)} \not\equiv 0 \indent\mbox{ (mod } p^2)
\end{equation*}
However,
\begin{eqnarray*}
p\alpha_1 &\equiv& 0 \indent\mbox{ (mod } p)\\
\Rightarrow \indent\indent F_{p\pi(p)} &\equiv& 0 \indent\mbox{ (mod } p^2)\indent\mbox{ (by Lemma 7)}\\
\Rightarrow \indent~~ F_{p\pi(p)+1} &\equiv& 1 \indent\mbox{ (mod } p^2)\indent\mbox{ (by Lemma 7)}
\end{eqnarray*}
Since $p\pi(p)$ is the least positive integer $g$ such that $F_g$ and $F_{g+1}$ are $0$ and $1$ modulo $p^2$ respectively, we conclude
$$
\pi(p^2) = p\pi(p)
$$

\emph{Case (ii):} \indent$\beta_1 \not\equiv 0 $ \indent(mod $p$)

In this case, $\beta_1$, $2\beta_1$, $3\beta_1$, ... $(p-1)\beta_1$ are not congruent to zero modulo $p$ as well.
Thus, by Lemma 7, we have for all positive integers $\xi<p$,
\begin{equation*}
F_{\xi\pi(p)+1} \not\equiv 1 \indent\mbox{ (mod } p^2)
\end{equation*}
However,
\begin{eqnarray*}
p\beta_1 &\equiv& 0 \indent\mbox{ (mod } p)\\
\Rightarrow \indent\indent F_{p\pi(p)} &\equiv& 0 \indent\mbox{ (mod } p^2)\indent\mbox{ (by Lemma 7)}\\
\Rightarrow \indent ~~F_{p\pi(p)+1} &\equiv& 1 \indent\mbox{ (mod } p^2)\indent\mbox{ (by Lemma 7)}
\end{eqnarray*}
Since $p\pi(p)$ is the least positive integer $g$ such that $F_g$ and $F_{g+1}$ are $0$ and $1$ modulo $p^2$ respectively, we conclude
$$
\pi(p^2) = p\pi(p)
$$

\emph{Case (iii):} \indent$\alpha_1 \equiv \beta_1 \equiv 0$ \indent(mod $p$)

In this case, we directly have by definition
\begin{eqnarray*}
F_{\pi(p)} &\equiv& 0 \indent\mbox{ (mod } p^2)\\
F_{\pi(p)+1} &\equiv& 1 \indent\mbox{ (mod } p^2)
\end{eqnarray*}
Hence, by again using the result $\pi(p)\mid\pi(p^2)$, we have $\pi(p^2) = \pi(p)$.

The Lemma is thus proved. $\blacksquare$

\section{A Few Equivalences}

Consider the following statements:
\\

\emph{\bf{Statement 2:}} For any given prime $p$, the highest power of $p$ dividing $F_{\kappa(p)}$ is 1.
\\

\emph{\bf{Statement 3:}} Let $m$ be a positive integer and $p$ be a prime. If $p^2 \mid F_m$, then $p$ is a proper divisor of $m$.
\\

\emph{\bf{Statement 4:}} For every prime $p$, we have,
$$
\pi(p^2) = p\pi(p)
$$
With the above statements in mind, we have the following Theorem:\\

\emph{\bf{Theorem:}} Statements 1, 2, 3 and 4 are equivalent.

\emph{\bf{Proof:}}  We prove the bidirectional implications: Statements 1 $\Leftrightarrow$ 2, 2 $\Leftrightarrow$ 3 and 2 $\Leftrightarrow$ 4.

To prove that Statement 1 implies Statement 2, we use the result in \cite{30},~\cite{31}:
  $$p\mid F_{p-\left(\frac{p}{5}\right)}$$
  According to Statement 1, $p^2 \nmid F_{p-\left(\frac{p}{5}\right)}$, which by Lemma 2 leads us to conclude $p^2\nmid F_{\kappa(p)}$.

To prove that Statement 2 implies Statement 1, we see that if the highest power of $p$ dividing $F_{\kappa(p)}$ is $1$ then, by Lemma 4, $\kappa(p^{2}) = p\kappa(p)$. But, $p\nmid~p-\left(\frac{p}{5}\right)$ except for the case $p=5$, which can be easily ruled out as 5 is not a Wall-Sun-Sun prime. Thus, we deduce $\kappa(p^{2})\nmid p-\left(\frac{p}{5}\right)$. Using Lemma 1, we conclude that $p^2\nmid F_{p-\left(\frac{p}{5}\right)}$ i.e. $p$ is not a Wall-Sun-Sun prime. Since, $p$ can be any arbitrary prime, we see that there does not exist any Wall-Sun-Sun prime.

To prove that Statement 2 implies Statement 3, we see that if $1$ is the highest power of $p$ dividing $F_{\kappa(p)}$, in accordance with Statement 2, then, by Lemma 4, $\kappa(p^{2}) = p\kappa(p)$. Now, by Lemma 1, if $p^2$ divides $F_m$ then $\kappa(p^2)$ i.e. $p\kappa(p)$ divides $m$. In other words, $p$ is a proper divisor of $m$.

To prove that Statement 3 implies Statement 2, we assume $p^2\mid F_{\kappa(p)}$. By Statement 3, $p$ is a proper divisor of $\kappa(p)$. So, we let $\kappa(p) = pc$ where $c > 1$. We have, by Lemma 3, $p\mid F_{p}$ or $p\mid F_{c}$. But both $p$ and $c$ are less than $\kappa(p)$, which according to definition is the least positive integer $n$ such that $p\mid F_n$. This is a contradiction. Hence $p^2\nmid F_{\kappa(p)}$. It follows that Statements 2 and 3 are equivalent.
\let\thefootnote\relax\footnotetext{$^\dagger$Though the equivalence between Statements 1 and 4 has been well established, we have provided an elementary proof.}

To prove that Statement 2 implies Statement 4$^\dagger$, we use Lemma 4 in conjunction with Statement 2 to note $\kappa(p^2) = p\kappa(p)$. But by Lemma 5, we know that:
\begin{equation*}
F_{p\kappa(p)+1} \equiv F_{\kappa(p)+1}^p \indent\mbox{ (mod } p^2)
\end{equation*}
We have, as a consequence:
\begin{equation*}
\Omega_{p^2}(F_{\kappa(p^2)+1}) = \Omega_{p^2}(F_{\kappa(p)+1}^p)
\end{equation*}
If we denote the quantity on either side as $\omega$, then
$$
F_{\kappa(p)+1}^{p\omega} \equiv 1 \indent\mbox{ (mod } p^2)\\
$$
Hence, we would also have
$$
F_{\kappa(p)+1}^{p\omega} \equiv 1 \indent\mbox{ (mod } p)\\
$$By Fermat's Little Theorem, we have
\begin{eqnarray*}
F_{\kappa(p)+1}^{\omega} &\equiv& F_{\kappa(p)+1}^{p\omega}\\
&\equiv& 1 \indent\mbox{ (mod } p)
\end{eqnarray*}
This means that $\Omega_{p}(F_{\kappa(p)+1})\mid \omega$ i.e. $\Omega_{p}(F_{\kappa(p)+1})\mid \Omega_{p^2}(F_{\kappa(p^2)+1})$ whence we get,
$$
\Omega_{p}(F_{\kappa(p)+1}) \le \Omega_{p^2}(F_{\kappa(p^2)+1})
$$
But we have,
\begin{eqnarray*}
\kappa(p^2) &=& p\kappa(p)\\
\Rightarrow\indent\indent\indent\indent~~~ \kappa(p^2) &>& \kappa(p)\\
\Rightarrow\indent\kappa(p^2)\Omega_{p^2}(F_{\kappa(p^2)+1}) &>& \kappa(p)\Omega_{p}(F_{\kappa(p)+1})
\end{eqnarray*}
And as a consequence of Lemma 6 we have,
$$
\pi(p) < \pi(p^2)
$$
Hence by Lemma 8,
$$
\pi(p^2) = p\pi(p)
$$
To prove that Statement 4 implies Statement 2, we merely note that if every Fibonacci number divisible by a given prime $p$ was also divisible by $p^2$, then $\pi(p^2)$ would be equal to $\pi(p)$. Hence, Statement 4 implies Statement 2.

The equivalence is thus proved. $\blacksquare$

\section{Heuristic Arguments}
Firstly, some exciting results have been proved by A. S. Elsenhans and J. Jahnel in \cite{power14} and we would request the readers to go through them.

An investigation regarding Wall-Sun-Sun Prime Conjecture carried out in \cite{power14} makes us believe that it \emph{might} be true. The popular version of the conjecture is its equivalent Statement 4. And from Lemma 8 we have that the conjecture implies that there are no solutions to the equation below in prime numbers:
$$
\pi(p^2)=\pi(p)
$$
However it would be interesting to find solutions to the above equation over all positive integers. Regarding this, we conjecture that the only solutions to the equation:
$$
\pi(n^2)=\pi(n),~~~ \forall n\in\mathbb{N}
$$
are $n=6$ and $n=12$. Although no clear reason presents itself to us now, as to why the number 6 has such an interesting relationship with its Pisano period, we can speculate why 12 follows it up. The Pisano period function, $\pi$ bears certain striking similarities to Euler's totient function, $\phi$. As indicated by computer investigation, for instance, both seem to obey similar relations:
$$\phi(p^n)=p^{n-1}\phi(p)\indent \cite{NT}$$
$$\pi(p^n)=p^{n-1}\pi(p)\indent \cite{power14}$$
for all primes $p$. Further results such as:
$$\begin{array}{lcl}
\phi(mn)=\phi(m)\phi(n), &\mbox{if }&\gcd(m,n)=1, \forall m,n\in\mathbb{N}\indent\cite{NT}\\
\pi(mn)=\mbox{lcm}(\pi(m),\pi(n)), &\mbox{if }&\gcd(m,n)=1, \forall m,n\in\mathbb{N}\indent\cite{power14}\\
\end{array}$$
confirm that there might be deeper links between the two functions. Now considering the above equations, it is easier to appreciate why $\pi(6)=\pi(12)=\pi(6^2)=\pi(12^2).$

On a different note, it has been intuitively argued that $\pi(p^2)=p\pi(p)$ for prime $p$ \cite{power14}. So it is reasonable to expect every prime $p$ to divide $\pi(p^2)$. However  for small values of $n$, it can be verified that:
$$n\mid\pi(n^2),~~~~~~n\in\mathbb{N} $$
Keeping in mind Lemma 4 and certain results mentioned in \cite{power14}, we claim that:
$$n\mid\pi(n^2),~~~~~~\forall n\in\mathbb{N} $$
If we see the above claim in the light of Lemma 6, we attain a better insight into the heart of the problem, which only becomes more compelling when we bound $\pi(n^2)$ by,
$$\pi(n)\le\pi(n^2)\le n\pi(n),~~~~~~\forall n\in\mathbb{N}$$ The proof of the inequality is omitted here, but we encourage the reader to prove them. (\emph{Hint}: Use Pigeonhole Principle). Note that no easily detectable pattern emerges, as to when the equality holds for the upper bound. Also, we have already conjectured regarding the condition when the equality holds for the lower bound. Now assuming the above bounds on $\pi(n^2)$, the claim that $ n\mid\pi(n^2), ~\forall n\in\mathbb{N}$, becomes even more intriguing. Summarizing, we conjecture the following statements:
\begin{enumerate}
\item The only solutions for the equation $\pi(n^2) = \pi(n)$ over positive integers are 6 and 12.
\item $n\mid\pi(n^2),~\forall n\in\mathbb{N}$
\end{enumerate}
\section{Acknowledgments}
We would like to thank Suryateja Gavva$^\dag$ for helping us verify the proofs, checking for any lapses of logic and providing us with the much needed motivation to complete this paper. We would also like to thank Vihang~Mehta$^\ddag$ for helping us cross-refer certain sources.
\let\thefootnote\relax\footnotetext{$\dag$ Suryateja Gavva is a second-year undergraduate at the Indian Institute of\\ Technology (IIT), Bombay.}
\let\thefootnote\relax\footnotetext{$\ddag$ Vihang Mehta is a third-year undergraduate at Brown University.}

\end{document}